\documentclass[11pt]{article}   

\usepackage{amsfonts}

\begin{document}

\parindent = 0pt \baselineskip = 22pt   \parskip = \the\baselineskip
\newcommand{\la}{\lambda} \newcommand{\Aa}{\mathbb{A}}
\newcommand{\RR}{\mathbb{R}} \newcommand{\FF}{\mathbb{F}}
\newcommand{\CC}{\mathbb{C}} \newcommand{\QQ}{\mathbb{Q}}
\newcommand{\HH}{\mathbb{H}} \newcommand{\ZZ}{\mathbb{Z}}
\newcommand{\NN}{\mathbb{N}}

\newcommand{\cB}{{\cal B}} \newcommand{\cP}{{\cal P}} \newcommand{\cV}{{\cal V}}
\newcommand{\cC}{{\cal C}} \newcommand{\cS}{{\cal S}} \newcommand{\cW}{{\cal W}}
\newcommand{\cF}{{\cal F}} \newcommand{\Real}{\mathop{\rm Re}}
\newcommand{\llim}{\mathop{\rm l{.}i{.}m{.}}}

\newtheorem{thm}{Theorem}[section] 
\newtheorem{thmdef}[thm]{Theorem and Definition} 
\newtheorem{prop}[thm]{Proposition} 
\newtheorem{lem}[thm]{Lemma}
\newtheorem{cor}[thm]{Corollary}

\newtheorem{defi}[thm]{Definition}

\newtheorem{pre-note}[thm]{Note}
\newenvironment{note}{\begin{pre-note}\rm}{\end{pre-note}}

\newenvironment{proof}{\bf Proof\ \rm}{$\;\bullet$}

\begin{titlepage}
\rightline{math.NT/0103058} 
\rightline{(submitted)} 
\indent
\begin{center}


{\bf\large A lower bound in an approximation problem involving the zeros of the
Riemann zeta function}\\

\vskip 1cm

Jean-Fran\c{c}ois Burnol\\

March 2001\\

\end{center}

\vskip 1cm

{\bf Abstract:} We slightly improve the lower bound of B\'aez-Duarte, Balazard,
Landreau and Saias in the Nyman-Beurling formulation of the Riemann Hypothesis
as an approximation problem. We construct Hilbert space vectors which could
prove useful in the context of the so-called ``Hilbert-P\'olya idea''.

\vfill

{\parskip = 0pt\parindent = 100bp\baselineskip=14bp
Author's affiliation:\par
Jean-Fran\c{c}ois Burnol\par
Universit\'e de Nice \--\ Sophia Antipolis\par
Laboratoire J.-A. Dieudonn\'e\par
Parc Valrose\par
F-06108 Nice Cedex 02\par
France\par 
electronic mail: burnol@math.unice.fr\par}

\end{titlepage}

\setcounter{page}{2}

\tableofcontents













\section{Introduction}

In \cite{Co1} and the subsequent paper \cite{Co2}, Connes gave a rather
intrinsic construction of a Hilbert space intimately associated with the zeros
of the Riemann zeta function on the critical line. But the zeros having
multiplicities higher than a certain level (which is a parameter in Connes's
construction),  have (if they at all exist) their contributions limited to that
level, and not to the extent given by their natural multiplicities. Thus
subsists the problem of a natural definition of a so-called ``Hilbert-P\'olya
space'', with orthonormal basis indexed by the zeros $\rho$ of $\zeta$ and
integers $k$ varying from $0$ to $m_\rho - 1$ where $m_\rho$  is the
multiplicity of $\rho$. We do not solve that problem here but we do propose a
rather natural construction of Hilbert space vectors $X^\la_{\rho, k}$,
$\zeta(\rho)=0$, $k<m_\rho$, which in the limit when the parameter $\la$ goes to
$0$ become perpendicular (when they correspond to distinct zeros. The vectors
corresponding to a multiple root $\rho$ are independent but need to be
orthogonalized.) As in Connes's constructions these vectors live in a quotient
space. Controlling the limit $\la\to0$ to obtain a so-called Hilbert-P\'olya
space probably involves considerations from mathematical scattering theory (we
have previously studied in \cite{Bu3}, \cite{Bu4} some connections with the
problems of $L-$functions.)

The context in which our construction takes place is that of the Nyman-Beurling
formulation of the Riemann Hypothesis as an approximation problem \cite{Nym},
\cite{Beu}. Let $K = L^2(]0,\infty[,dt)$ (over the complex numbers), let $\chi$
be the indicator function of the interval $]0,1]$, and let $\rho$ be the
function ``fractional part'' (the letter $\rho$ is also used to refer to a zero
of the Riemann zeta function, hopefully no confusion will arise.) Let
$0<\lambda<1$ and let $\cB_\lambda$ be the sub-vector space of $K$ consisting of
the finite linear combinations of the functions $t\mapsto \rho({\theta\over
t})$, for $\lambda\leq\theta\leq1$.

\begin{thm}[Nyman \cite{Nym}, Beurling \cite{Beu}]
The Riemann Hypothesis holds if and only if $$\chi\in
\overline{\bigcup_{0<\lambda<1} \cB_\lambda}$$
\end{thm}

Actually we are following \cite{Bal} here in using a slight variant of the
original Nyman-Beurling formulation. It is a disappointing fact that this
theorem can be proven without leading to any new information whatsoever on the
zeros lying on the critical line (basically what is at works is the
factorization of functions belonging to the Hardy space of a half-plane
\cite{Hof}.) The following is thus rather remarkable:

\begin{thm}[B\'aez-Duarte, Balazard,
Landreau and Saias \cite{Bal}] Let us write $D(\la)$ for the Hilbert-space
distance $\inf_{f\in\cB_\lambda} \| \chi -f \|$.   We have
$$\liminf_{\la\to0} D(\la)\sqrt{\log({1\over\la})} \geq \sqrt{\sum_\rho {1\over
|\rho|^2}}$$
\end{thm}

If the Riemann Hypothesis fails this result is true but trivial as
the left-hand side then takes the value $+\infty$. So we will assume that the
Riemann Hypothesis holds. The sum on the right-hand side is over all non-trivial
zeros $\rho$ of the zeta function, counted \emph{only once} independently of
their multiplicities $m_\rho$. We prove the following:

\begin{thm} We have:
$$\liminf_{\la\to0} D(\la)\sqrt{\log({1\over\la})} \geq \sqrt{\sum_\rho
{m_\rho^2\over |\rho|^2}}$$
\end{thm}

So the zeros are counted according to the \emph{square} of their
multiplicities. To prove this lower bound we will construct remarkable Hilbert
space vectors $X^\la_{\rho, k}\;$, $\zeta(\rho)=0$, $k<m_\rho$ and use them to
control $D(\la)$.  The following ``toy-model'' gives us reasons to expect that
the lower bound in fact gives the exact order of decrease of $D(\la)$:

\begin{thm}\label{theotoy}
Let $Q(z) = \prod_\alpha (1-\overline{\alpha}\cdot z)^{m_\alpha}$ be a
polynomial of degree $q\geq1$ will all its roots $\alpha$ on the unit circle
(the root $\alpha$ having multiplicity $m_\alpha$). Let $P(z)$ be an arbitrary
polynomial. Let
$$E(N,P) := \inf_{\deg(A)\leq N} \int_{S^1} |P(z) -
Q(z)A(z)|^2{d\theta\over2\pi}$$ We have as $N$ goes to infinity:
$$\lim\ N\,E(N,P) = \sum_\alpha {m_\alpha^2}\;|P(\alpha)|^2$$
\end{thm}

\section{The prediction error for a singular MA(q)}

As motivation for our result we first consider a simpler approximation
problem, in the context of the Hardy space of the unit disc rather than the
Hardy space of a half-plane. Let $Q(z) = \prod_\alpha (1-\overline{\alpha}\cdot
z)^{m_\alpha}$ be a polynomial of degree $q\geq1$ will all its roots $\alpha$ on
the unit circle (the root $\alpha$ having multiplicity $m_\alpha$ so that $q =
\sum_\alpha m_\alpha$.) Let us define:
$$E(N) := \inf_{\deg(A)\leq N} \int_{S^1} |1 - Q(z)A(z)|^2{d\theta\over2\pi}$$
The measure ${d\theta\over2\pi}$ is the rotation invariant probability measure
on the circle $S^1$, with $z =\exp(i\theta)$. The minimum is taken over all
complex polynomials $A(z)$ with degree at most $N$. We are guaranteed that
$\lim_{N\to\infty} E(N) = 0$ as $Q(z)$ is an outer factor (\cite{Hof}). More
precisely:

\begin{thm} As $N$ goes to infinity we have:
$$\lim\ N\,E(N) = \sum_\alpha {m_\alpha^2}$$
\end{thm}

\begin{note} 
In case $Q(z)$ has a root in the open unit disc then $E(N)$ is bounded below by
a positive constant. In case $Q(z)$ has all its roots outside the open unit
disc, then the result above holds but only the roots on the unit circle
contribute. Finally if all its roots are outside the closed unit disc then the
decrease is exponential: $E(N) = O(c^N)$, with $c<1$.
\end{note}

The theorem, although not stated explicitely there, is easily extracted from the
work of Grenander and Rosenblatt \cite{Gre}. They state an $O({1\over N})$
result, in a much wider set-up than the one considered here (which is limited to
simple-minded $q$-th order moving averages.) Unfortunately the $O({1\over N})$
bound is now believed not to be systematically true under their hypotheses (as
is explained in \cite{Nev}; I thank Professor W.~Van~Assche for pointing out
this fact to me.) Nevertheless their technique of proof goes through smoothly in
the case at hand and yields the exact asymptotic result as stated above. We only
sketch briefly the idea, as nothing  beyond the tools used in \cite{Gre} is
needed.

We point out in passing that it is of course possible to express $E(N)$
explicitely in terms of the Toeplitz determinants for the measure $d\mu =
|Q(\exp(i\theta))|^2 {d\theta\over2\pi}$. But already for an $MA(2)$ this gives
rise to unwieldy computations\dots. Rather: let $\cP_{N}$ be the vector space of
polynomials of degrees at most $N+q$, let $\cV_{N}$ be the subspace of
polynomials divisible by $Q(z)$, and let $\cW_N$ be its $q$-dimensional
orthogonal complement. Then $E(N)$ is the squared norm of the orthogonal
projection of the constant function $1$ to $\cW_N$. A spanning set in $\cW_N$ is
readily identified: to each root $\alpha$ one associates $Y_{\alpha,0}^N$,
$Y_{\alpha,1}^N$, \dots, $Y_{\alpha,m_\alpha -1}^N$ defined as
$$Y_{\alpha,0}^N := 1 + \overline{\alpha} z + \dots +\overline{\alpha}^{N+q}
z^{N+q}$$
$$Y_{\alpha,1}^N := z + 2\overline{\alpha} z^2 + \dots +
(N+q)\overline{\alpha}^{N+q-1} z^{N+q}$$ and similarly for $k=2, \dots, m_\alpha
- 1$. We can then express $E(N)$ using a Gram formula in terms of (the inverse)
of the positive matrix (of fixed size $q\times q$ but depending on $N$) built
with the scalar products of the $Y$'s. It turns out that in the limit when $N$
goes to infinity and after the rescaling $Y_{\alpha,k}^N\mapsto X_{\alpha,k}^N
:= N^{-k-1/2}Y_{\alpha,k}^N$ the Gram matrix decomposes into Cauchy blocks
$(1/(i+j+1))_{0\leq i,j < m_\alpha}$ of size $m_\alpha$, one for each root
$\alpha$. It is known from Cauchy that the top-left element of the inverse
matrix is $m_\alpha^2$. This is how ${\sum_\alpha {m_\alpha^2}\over N}$ arises,
after keeping track of the scalar products $(1, X_{\alpha,k}^N)$. Instead of the
constant polynomial $1$ we could have looked at the approximation rate to an
arbitrary polynomial $P(z)$. The proof just sketched applies identically and
gives the Theorem \ref{theotoy} from the Introduction.

\section{Invariant analysis and a construction of B\'aez-Duarte}

The Mellin transform $f(t)\mapsto \widehat{f}(s)=\int_{t>0} f(t)t^{s-1}dt$
establishes the Plancherel isometry between $K = L^2(]0,\infty[,dt)$ and
$L^2(s={1\over2}+i\tau,{d\tau\over2\pi})$, with inverse $F(s)\mapsto
\int_{s=1/2+i\tau} F(s) t^{-s}{d\tau\over2\pi}$. Let $a(s)$ be a measurable
function of $s$ (as a rule when using the letter $s$ we implicitely assume
$\Real(s) = {1\over2}$. We will use letters $w$ and $z$ for general complex
numbers.) If $a(s)$ is essentially bounded then $F(s)\mapsto a(s)F(s)$ defines a
bounded operator on $K$ which commutes with the unitary group
$D_{\theta}:f(t)\mapsto {1\over\sqrt{\theta}}f({t\over{\theta}})$, and all
bounded operators commuting with the $D_\theta$ ($0<\theta<\infty$) are obtained
in such a manner. More generally all \emph{closed} invariant operators are
associated to a measurable multiplier $a(s)$ (finite almost everywhere, but not
necessarily essentially bounded). For the details of this technical statement,
see \cite{Bu5}.

For example the Hardy averaging operator $M:f(t)\mapsto {1\over t}\int_{]0,t]}
f(u)du$ corresponds to the spectral multiplier $1\over 1-s$. The operator $1 -
M$ corresponds to the spectral multiplier $s\over s-1$ and is thus
unitary. Another (see \cite{Bu2}) remarkable invariant operator is the (even)
``Gamma'' operator $\Gamma_+ = \cF_+ I$. Here $I$ is the inversion $f(t)\mapsto
{1\over t}f({1\over t})$ and $\cF_+$ is the additive Fourier transform as
applied to even functions (the cosine transform). The multiplier associated to
$\Gamma_+$ is the (Tate) function
$$\gamma_+(s) = \pi^{{1\over2}-s}{\Gamma(s/2)\over\Gamma((1-s)/2)} =
{\zeta(1-s)\over\zeta(s)} = 2^{1-s}\pi^{-s}\cos({\pi\,s\over2})\Gamma(s) =
(1-s)\int_0^\infty u^{s-1}\,{\sin(2\pi u)\over\pi u}\,du$$
A further invariant operator is the operator $U$
introduced by B\'aez-Duarte \cite{Bae} in connection with the Nyman-Beurling
formulation of the Riemann Hypothesis: its spectral multiplier is ${s\over
1-s}{\zeta(1-s)\over\zeta(s)}$, so $U = (M-1)\cF_+ I = \cF_+ I (M-1)$.

From the results recalled above on invariant operators, we see that invariant
orthogonal projectors correspond to indicator functions of measurable sets on
the critical line. So a function $f(t)$ is such that its multiplicative
translates $D_{\theta}(f)$ ($0<\theta<\infty$) span $K$ if and only if
$F(s)=\widehat{f}(s)$ is almost everywhere non-vanishing (Wiener's
$L^2$-Tauberian Theorem.) In that case the phase function
$$U_f(s) = {\ \overline{F(s)}\ \over F(s)}$$ is almost everywhere defined and of
modulus $1$. It thus corresponds to an invariant unitary operator, also denoted
$U_f$.

Let us introduce the \emph{anti-unitary} ``time-reversal'' operator $J$  acting
on $K$ as $g\mapsto \overline{I(g)}$. The operator $U_f$ commutes with the
contractions-dilations, is unitary, and sends $f$ to $J(f)$. We call this the
\emph{B\'aez-Duarte construction} as it appears in \cite{Bae} (up to some
non-essential differences) in relation with the Nyman-Beurling problem (the
phase function arises in other contexts, especially in scattering theory.)

To relate this with the operator $U = (M-1)\cF_+ I$, one needs the formula
$${\zeta(s)\over s} = - \int_0^\infty \rho({1\over t}) t^{s-1} dt$$ which is
fundamental in the Nyman-Beurling context. This formula shows that $U$ is the
phase operator associated with $\rho({1\over t})$.

Generally speaking, the
operators $U_f$ are related to the Hardy spaces $\HH^2 = L^2(]0,1],dt)$ and
${\HH^2}^\perp = L^2([1,\infty[,dt)$ (we will also use the notation $\HH^2$ for the
Mellin transform of $L^2(]0,1],dt)$.) Indeed the time-reversal $J$ is an
isometry (anti-unitary) between $\HH^2$ and ${\HH^2}^\perp$. Let us assume that
the function $f$ belongs to $\HH^2$. The operator $U_f$ has the same effect as
$J$ on $f$, but contrarily to $J$ is an \emph{invariant} operator. This puts the
space $\cB_\la(f)$ (of finite linear combinations of contractions $D_\theta(f)$
for $\la\leq\theta\leq1$) isometrically in a new light as a subspace of
$L^2([\la,\infty[,dt)$. The marvelous thing is that in this new incarnation it
appears to be sometimes possible to find vectors
orthogonal to $\cB_\la(f)$ and thus to get some control on $\cB_\la(f)$ as $\la$
decreases (as in the Grenander-Rosenblatt method.)

\section{The vectors $Y^\la_{s,k}$}

To get started on this we first replace the $L^2$ function $-{\zeta(s)\over s}$
with an element of $\HH^2$. This is elementary:

\begin{prop}[\cite{Bu4}, \cite{Ehm}]
The function $Z(s) = {s-1\over s}{\zeta(s)\over s}$ belongs to $\HH^2$. Its
inverse Mellin transform $A(t)$ is given by the formula
$$A(t) = [{1\over t}]\log(t) + \log([{1\over t}]\,!) + [{1\over t}]$$ One has
(\/for $0<t\leq 1$) $A(t)= {1\over2}\log({1\over t}) + O(1)$.
\end{prop}

The B\'aez-Duarte construction will then associate to $A(t)$ the operator $V$ with
spectral multiplier
$$V(s) = \left({s\over 1-s}\right)^3\;{\zeta(1-s)\over \zeta(s)}$$ 
so that
$$V = (1 - M)^2 \cdot U$$ 
This last representation will prove useful as it allows to use the formulae
related to $U$ from \cite{Bae} and \cite{Bal}. Let $\cC_\la$ ($0<\la<1$) be
the sub-vector space of $\HH^2$ of linear combinations of the contractions
$D_\theta(A)$ for $\la\leq\theta\leq1$. The function ${s-1\over s}{1\over s} =
{1\over s} - {1\over s^2}$ is the Mellin transform of
$\chi_1(t):=(1+\log(t))\chi(t)$. The quantity $D(\lambda)$ considered by
B\'aez-Duarte, Balazard, Landreau and Saias is thus the Hilbert space distance
between $\chi_1(t)$ and $\cC_\la$. To bound it from below we will exhibit
remarkable Hilbert space vectors $X^\lambda_{\rho,k}$ indexed by the zeros of
the Riemann zeta function and perpendicular to $\cC_\la$. We then compute the
exact asymptotics of the orthogonal projection of $\chi_1$ to the vector spaces
spanned by the $X^\lambda_{\rho,k}$, for a finite set of roots, exactly as in
the Grenander-Rosenblatt method.

To each complex number $w$ and natural integer $k\geq0$ we associate the funtion
$\psi_{w,k}(t) = (\log({1\over t}))^k\, t^{-w}\,\chi(t)$ on $]0,\infty[$. For
$\Real(w)<1$ it is integrable, for $\Real(w)<{1\over2}$ it is in $K$. Let
$Q_\la$ be the orthogonal projector from $K$ onto $L^2([\la,\infty[)$. The main
point of this paper is the following:

\begin{thmdef}
For each $0<\la\leq 1$, each $s$ on the critical line, and each integer $k\geq0$
the $L^2$-limit in $K$ of $V^{-1}Q_\la V(\psi_{w,k})$ exists as $w$ tends to $s$
from the left half-plane:
$$Y^\la_{s,k}:= \llim_{w\to s\atop\Real(w)<{1\over2}}\ V^{-1}Q_\la
V(\psi_{w,k})$$ For each $\la\leq\theta\leq1$ the scalar products between
$D_\theta(A)$ and the vectors $Y^\la_{s,k}$ are:
$$\la\leq\theta\leq1\ \Rightarrow\ (D_\theta(A), Y^\la_{s,k}) = \left(-{d\over
ds}\right)^k\;\theta^{s - {1\over 2}} Z(s)$$
\end{thmdef}

\begin{note}
The proof shows the existence of an analytic continuation in $w$ accross the
critical line, but we shall not make use of this fact.
\end{note}

Clearly one has the following statement as an immediate consequence:

\begin{cor}
Let $0<\la<1$. The vector $Y^\la_{s,k}$ is perpendicular to $\cC_\la$ if and
only if $\zeta^{(j)}(s) = 0$ for all $j\leq k$, if and only if $s$ is a zero
$\rho$ of the zeta function and $k<m_\rho$.
\end{cor}

\begin{note}
Our scalar products $(f,g)$ are complex linear in the first factor and
conjugate-linear in the second factor.
\end{note}

\begin{note}
The operator ${d\over ds}$ when applied to a not necessarily analytic function on
the critical line is defined to act as ${1\over i}{d\over d\tau}$ (where $s =
{1\over2} + i\tau$.)
\end{note}

\begin{proof}
The proof of existence will be given later. Here we check the statement
involving the scalar product, assuming existence. The following holds for
$\la\leq\theta\leq1$ and $\Real(w)<{1\over2}$:
\begin{eqnarray*}
(V^{-1}Q_\la V(\psi_{w,k}), D_\theta(A)) &=& (Q_\la V(\psi_{w,k}),
VD_\theta(A))\cr &=& (Q_\la V(\psi_{w,k}), D_\theta\cdot V(A))\cr &=&
(V(\psi_{w,k}), Q_\la\cdot D_\theta\cdot J(A))\cr &=& (V(\psi_{w,k}),
D_\theta\cdot V(A))\cr &=& (\psi_{w,k}\;, D_\theta(A))\cr &=& \left({d\over
dw}\right)^k\ (\psi_{w,0}\;, D_\theta(A))\cr &=& \left({d\over dw}\right)^k\
(D_{\theta}^{-1}(t^{-w}\,\chi(t)), A)\cr &=& \left({d\over dw}\right)^k\
(\theta^{1/2 - w}t^{-w}\,\chi(\theta t), A)\cr &=& \left({d\over dw}\right)^k\
\theta^{1/2 - w} \int_{]0,1]} t^{-w} \overline{A(t)}\,dt \cr
\end{eqnarray*}
Taking the limit when $w\to s$ gives
\begin{eqnarray*}
(Y^\la_{s,k}, D_\theta(A)) &=& \left({d\over ds}\right)^k\ \theta^{1/2 - s}
\int_{]0,1]} t^{-s} \overline{A(t)}\,dt\cr &=& \left({1\over i}{d\over
d\tau}\right)^k\ \theta^{-i\tau}\int_{]0,1]}t^{-{1\over2} -
i\tau}\;\overline{A(t)}\,dt
\end{eqnarray*}
Taking the complex conjugate:
\begin{eqnarray*}
(D_\theta(A), Y^\la_{s,k}) &=&  \left({i}{d\over d\tau}\right)^k\
\theta^{i\tau}\int_{]0,1]}t^{-{1\over2} + i\tau}\;A(t)\,dt\cr &=& \left(-{d\over
ds}\right)^k\;\theta^{s - {1\over 2}} \int_{]0,1]}t^{s-1}\;A(t)\,dt\cr &=&
\left(-{d\over ds}\right)^k\;\theta^{s - {1\over 2}} Z(s)
\end{eqnarray*}
which completes the proof (assuming existence.)
\end{proof}

To prove the existence we will use in an essential manner the key {\bf Lemme 6}
from \cite{Bal}. We have seen that $V = (1 - M)^2 U$ where $M$ is the Hardy
averaging operator and $U$ the B\'aez-Duarte operator. The spectral function
$U(s)$ extends to an analytic function $U(w)$ in the strip $0<\Real(w)<1$. We
need pointwise expressions for $V(\psi_{w,k})(t)$, $t>0$ (at first only
$\Real(w)<{1\over2}$ is allowed here). Thanks to the general study of $U$ given
in \cite{Bae}, we know that for $\Real(w)<{1\over2}$ the vector $U(\psi_{w,k})$
in $K$ is given as the following limit in square mean:
$$\llim_{\delta\to0}\int_\delta^1 (\log({1\over v}))^k\, v^{-w}\,{d\over
dv}{\sin(2\pi t/v)\over \pi t/v}\,dv$$ Following \cite{Bal}, with a slight
change of notation, we now study for each complex number $w$ with $\Real(w)<1$
(and each integer $k\geq0$) the \emph{pointwise} limit as a function of $t>0$
for $\delta\to 0$:
$$\varphi_{w,k}(t) := \lim_{\delta\to0}\int_\delta^1 (\log({1\over v}))^k\,
v^{-w}\,{d\over dv}{\sin(2\pi t/v)\over \pi t/v}\,dv$$

\begin{thm}[\cite{Bal}]\label{theo2}
Let $k=0$. For each $t>0$ and $\Real(w)<1$ the pointwise limit defining $\varphi_{w,0}(t)$
exists. It is holomorphic in $w$ for each fixed $t$. When $w$ is restricted to a
compact set in $\Real(w)<1$, one has uniformly in $w$ the bound
$\varphi_{w,0}(t) = O({1\over t})$ on $[1,\infty[$. Uniformly with respect to
$w$ satisfying $0<\Real(w)<1$ one has $ \varphi_{w,0}(t)  = U(w)\,t^{-w} + O(1)$
on $0<t\leq 1$.
\end{thm}

\begin{proof}
Everything is either stated explicitely in \cite{Bal}, Lemme 6 and Lemme 4, or
follows from their proofs. We will give more details for $k\geq1$ as this is not
treated in \cite{Bal}.
\end{proof}

\begin{cor}
For each $w$ in the critical strip $0<\Real(w)<1$ the Hardy operator $M: f(t)\to
{1\over t}\int_0^t f(v)\,dv$ can be applied arbitrarily many times to
$\varphi_{w,0}(t)$. The functions $M^L(\varphi_{w,0})$ ($L\in\NN$) are
$O({(1+\log(t))^L\over t})$ on $[1,\infty[$, uniformly with respect to $w$ when
it is restricted to a compact subset of the open strip, and satisfy on
$t\in\,]0,1]$ the estimate $M^L(\varphi_{w,0})(t) = \left({1\over 1 -
w}\right)^L U(w)\,t^{-w} + O(1)$, uniformly with respect to $w$.
\end{cor}

\begin{proof}
A simple recurrence.
\end{proof}

We thus obtain:

\begin{cor}
The vectors $Y^\la_{s,0}$ exist (for $\Real(s) = {1\over2}$). One has the
estimates:
$$V(Y^\la_{s,0})(t) = O({(1+\log(t))^2\over t})\qquad (t\in [1,\infty[)$$
$$V(Y^\la_{s,0})(t) = V(s)\;t^{-s} + O(1)\qquad (\lambda<t\leq 1)$$
$$V(Y^\la_{s,0})(t) = 0\qquad (0<t<\lambda)$$ uniformly with respect to $s$ when
its imaginary part is bounded.
\end{cor}

\begin{thm}\label{pretheo}
Let $k\geq 1$. For each $t>0$ and $\Real(w)<1$ the pointwise limit defining
$\varphi_{w,k}(t)$ exists. It is holomorphic in $w$ for each fixed $t$. When $w$
is restricted to a compact set in $\Real(w)<1$, one has uniformly in $w$ the
bound $\varphi_{w,k}(t) = O({1\over t})$ on $[1,\infty[$. Uniformly for
$0<\Real(w)<1$ one has $\varphi_{w,k}(t)  = \left({d\over dw}\right)^k
(U(w)\,t^{-w}) + O(1)$ on $0<t\leq 1$.
\end{thm}

\begin{proof}
The formula defining $\varphi_{w,k}(t)$ is equivalent to (after integration by
parts and the change of variable $u=1/v$):
$$\varphi_{w,k}(t) = \lim_{\Lambda\to\infty} {1\over\pi\,t}\int_1^\Lambda
(k+w\log(u))\left(\log(u)\right)^{k-1}\,u^{w-1}\sin(2\pi t\,u)\,{du\over u}$$
This proves the existence of $\varphi_{w,k}(t)$, its analytic character in $w$,
and the uniform $O({1\over t})$ bound on $[1,\infty[$. The formula can be
rewritten as:
$$\varphi_{w,k}(t) = \left({d\over dw}\right)^k {w\,\over\pi\,t}\int_1^\infty
u^{w-1}\sin(2\pi t\,u)\,{du\over u}$$ When $w$ is in the critical strip the
integral $\int_0^\infty u^{w-1}\sin(2\pi t\,u)\,{du\over u}$ is absolutely
convergent and its value is $t^{1-w}\int_0^\infty u^{w-1}\sin(2\pi u)\,{du\over
u} = {1\over 1-w}(2\pi\,t)^{1-w}\cos({\pi w\over2})\Gamma(w)$ from well-known
integral formulae, so that:
$$\varphi_{w,k}(t) = \left({d\over dw}\right)^k \left({w\over
1-w}2^{1-w}\pi^{-w}\cos({\pi w\over2})\Gamma(w) t^{-w} -
{w\,\over\pi\,t}\int_0^1 u^{w-1}\sin(2\pi t\,u)\,{du\over u}\right)$$ The first
term is $\left({d\over dw}\right)^k (U(w)\,t^{-w})$ and the second term can be
explicitely evaluated using the series expansion of $\sin(2\pi t\,u)$ with the
final result
$$\varphi_{w,k}(t) = \left({d\over dw}\right)^k (U(w)\,t^{-w}) + 2 (-1)^k \;
k!\sum_{j\geq1}(-1)^{j}{(2\pi t)^{2j}\over (2j+1)!}{2j\over (w+2j)^{k+1}}$$
which shows $\varphi_{w,k}(t)  = \left({d\over dw}\right)^k (U(w)\,t^{-w}) +
O(1)$, on $0<t\leq 1$, uniformly for $0<\Real(w)<1$.
\end{proof}

As was the case for $k=0$ we then deduce that the Hardy operator can be applied
arbitrarily many times to $\varphi_{w,k}$ for $0<\Real(w)<1$. The existence
of the $Y^\la_{s,k}$ follows.

\begin{thm}\label{theo}
Let $k\geq 0$. The vectors $Y^\la_{s,k}$ exist (for $\Real(s) = {1\over2}$). One
has the estimates:
\begin{eqnarray*}
V(Y^\la_{s,k})(t) &=& O({(1+\log(t))^2\over t})\qquad (t\in [1,\infty[)\cr
V(Y^\la_{s,k})(t) &=& \left({d\over ds}\right)^k (V(s)\;t^{-s}) + O(1)\qquad
(\lambda<t\leq 1)\cr V(Y^\la_{s,k})(t) &=& 0\qquad (0<t<\lambda)\cr
\end{eqnarray*}
the implied constants are independent of $\la$ and are uniform with respect to $s$
when its imaginary part is bounded.
\end{thm}

\begin{proof}
Clearly a corollary to \ref{pretheo}.
\end{proof}

\section{The vectors $X^\la_{\rho,k}$ and completion of the proof}

\begin{defi}
Let $0<\la<1$. To each zero $\rho$ of the Riemann zeta function on the critical
line, of multiplicity $m_\rho$, and each integer $0\leq k<m_\rho$ we associate
the Hilbert space vector
$$X^\la_{\rho,k} := \left(\log({1\over\lambda})\right)^{-{1\over2}-k}\cdot
Y^\la_{\rho,k}$$ where $Y^\la_{\rho,k} = \llim_{w\to s} V^{-1}Q_\la
V(\psi_{w,k})$, $V$ is the unitary operator $(M-1)^3 \cF_+\,I$,
$Q_\la$ is orthogonal projection to $L^2([\la,\infty[,dt)$, and $\psi_{w,k}(t) =
(\log({1\over t}))^k\, t^{-w}\;\chi(t)$.
\end{defi}

\begin{note}
Of course there is no reason except psychological to allow only zeros of the
Riemann zeta function at this stage.
\end{note}

\begin{thm}\label{theogram}
As $\la$ decreases to $0$ one has:
\begin{eqnarray*}
\lim_{\la\to0}\ (X^\la_{\rho_1,k}, X^\la_{\rho_2,l}) &=&
0\qquad(\rho_1\neq\rho_2)\cr \lim_{\la\to0}\ (X^\la_{\rho,k}, X^\la_{\rho,l})
&=& {1\over k+l+1}
\end{eqnarray*}
\end{thm}

\begin{proof}
To establish this we first consider, for $\Real(s_1) = \Real(s_2) = {1\over2}$:
$$\int_\la^1 \left(\log({1\over t})\right)^{j_1} t^{-s_1}\;\left(\log({1\over
t})\right)^{j_2} {t^{-(1-s_2)}}\,dt$$ If $s_1\neq s_2$ an integration by parts
shows that it is $O\left(\log({1\over\lambda})\right)^{j_1+j_2} $. On the other
hand when $s_1 = s_2$  its exact value is ${1\over j_1 + j_2 +
1}\left(\log({1\over\lambda})\right)^{j_1+j_2 + 1}$. With this information the
theorem follows directly from \ref{theo} as (for example) the leading divergent
contribution as $\la\to0$ to  $\left(V(Y^\la_{s,k}),\;V(Y^\la_{s,l})\right)$ is
$V(s)\overline{V(s)}\int_\la^1  \left({d\over ds}\right)^k\;t^{-s}\
\overline{\left({d\over ds}\right)^l\;t^{-s}}\,dt$ which gives ${1\over k + l +
1}\left(\log({1\over\lambda})\right)^{k+l+ 1}$. The rescaling $Y\mapsto X$ is
chosen so that a finite limit for $(X^\la_{\rho,k}, X^\la_{\rho,l})$ is
obtained. As the scalar products involving distinct zeros have a smaller
divergency, the rescaling let them converge to $0$.
\end{proof}

\begin{thm}\label{theoscal}
Let $\chi_1(t) = (1 + \log(t))\chi(t)$. As $\la$ decreases to $0$ one has:
\begin{eqnarray*}
\lim_{\la\to0}\ \sqrt{\log({1\over\la})}\;(\chi_1, X^\la_{\rho,k}) &=&
0\qquad(k\geq1)\cr \lim_{\la\to0}\ \sqrt{\log({1\over\la})}\;(\chi_1,
X^\la_{\rho,0}) &=& {\rho - 1\over\rho^2}\cr
\end{eqnarray*}
\end{thm}

\begin{proof}
We have $(1 - M)\chi_1 = \chi$, and $V = (1-M)^2\,U$ so
$V\chi_1 = (1-M)\,U\chi$. From \cite{Bae} we know that $U\chi$ is ${\sin(2\pi
t)\over \pi t}$ so $V\chi_1$ is the function ${\sin(2\pi t)\over \pi t} -
{1\over t}\int_0^t {\sin(2\pi v)\over \pi v}\,dv$. It is thus $0(t^2)$ as
$t\to0$, and from \ref{theo} we then deduce that the scalar products $(\chi_1,
Y^\la_{\rho,k})$ admit finite limits as $\la\to0$. This settles the case
$k\geq1$. For $k=0$, one uses the uniformity with respect to $w$ in \ref{theo2}
to get $$\lim_{\la\to0} (\chi_1, Y^\la_{\rho,0}) = \lim_{w\to\rho} (\chi_1,
\varphi_{w,0})$$ which gives $\lim_{w\to\rho} \int_0^1 (1+\log(t))\,t^{w-1}\,dt
= {1\over \rho} - {1\over\rho^2} = {\rho - 1\over\rho^2}$.
\end{proof}

We can now conclude the proof of our estimate.

\begin{thm}\label{mytheo}
We have:
$$\liminf_{\la\to0} D(\la)\sqrt{\log({1\over\la})} \geq \sqrt{\sum_\rho
{m_\rho^2\over |\rho|^2}}$$
\end{thm}

\begin{proof}
Let $R$ be a non-empty finite set of zeros. We showed that $D(\la)$ is the
Hilbert space distance from $\chi_1$ to $\cC_\lambda$, and that the vectors
$X^\la_{\rho,k}$ for $0\leq k<m_\rho$ are perpendicular to $\cC_\lambda$. So
$D(\la)$ is bounded below by the norm of the orthogonal projection of $\chi_1$
to the finite-dimensional vector space $H_R$ spanned by the vectors
$X^\la_{\rho,k}$, $0\leq k<m_\rho$, $\rho\in R$. This is given by a well-known
formula involving the inverse of the Gram matrix of the $X^\la_{\rho,k}$'s as
well as the scalar products $(\chi_1, X^\la_{\rho,k})$. From \ref{theogram} the
Gram matrix converges to diagonal blocks, one for each zero, given by Cauchy
matrices of sizes $m_\rho\times m_\rho$. From Cauchy we know that the top-left
element of the inverse matrix is $m_\rho^2$. Combining this with the scalar
products evaluated in \ref{theoscal} we get that the squared norm of the
orthogonal projection of $\chi_1$ to $H_R$ is asymptotically equivalent as
$\la\to0$ to $\sum_{\rho\in R} {m_\rho^2\over |\rho|^2}\over
\log({1\over\la})$. The proof is complete.
\end{proof}

We can apply our strategy to a fully singular MA(q) on the unit circle. The
relevant B\'aez-Duarte phase operator will then be (up to a non-important
constant of modulus 1) the operator of multiplication by $z^{-q}$ and it is
apparent that this leads to a proof equivalent to the one we gave in our
previous discussion, inspired by \cite{Gre}. In the case of the Nyman-Beurling
approximation problem for the zeta funtion, we expect in the quotient of $\HH^2$
by $\overline{\cC_\la}$ a ``continuous spectrum'' additionally to the ``discrete
spectrum'' provided by the (projection to $\HH^2$ of the) $X^\la_{\rho,k}$'s,
$\zeta(\rho)=0$, $k<m_\rho$. It is  tempting to speculate that the rescaling
will kill this continuous part as $\la\to0$, so that in the end only subsists a
so-called ``Hilbert-P\'olya'' space. This would appear to require \ref{mytheo}
to give the exact order of decrease of the quantity $D(\la)$ and the numerical
explorations reported by B\'aez-Duarte, Balazard, Landreau and Saias in
\cite{Bal} seem to support this.

\section{Acknowledgements}

I had been looking for the vectors $X^\lambda_{\rho,k}$, and with the {\bf Lemme
6} of \cite{Bal} they were suddenly there. I thank Michel Balazard for giving me
copies of \cite{Bae} and \cite{Bal} in preprint form.

%

\clearpage

\clearpage
\baselineskip=14pt\parskip=12pt

\begin{thebibliography}{99}

\bibitem{Bae} L.~B\'aez-Duarte,  \emph{A class of invariant unitary
operators}, Adv. in Math. {\bf 144} (1999), 1-12.

\bibitem{Bal} L.~B\'aez-Duarte, M.~Balazard, B.~Landreau and E.~Saias,
\emph{Notes sur la fonction $ \zeta $ de Riemann 3},  Adv. in Math.
{\bf 149} (2000), 130-144.

\bibitem{Beu} A.~Beurling, \emph{A closure problem related to the Riemann
Zeta--function\/}, Proc.\ Nat.\ Acad.\ Sci.\ {\bf 41} (1955), 312-314.


\bibitem{Bu2} J.-F.~Burnol, \emph{``Sur les formules explicites I: analyse
invariante''}, C. R. Acad. Sci. Paris {\bf 331} (2000), S\'erie I, 423-428.

\bibitem{Bu3} J.-F.~Burnol, \emph{Scattering on the p-adic field and a trace
formula}, International Mathematical Research Notices 2000:{\bf 2} (2000), 57-70.

\bibitem{Bu4} J.-F.~Burnol, \emph{An adelic causality problem related to
abelian $L-$functions}, Journal of Number Theory. To appear.

\bibitem{Bu5} J.-F.~Burnol, \emph{Quaternionic gamma functions and their
logarithmic derivatives as spectral functions}, Mathematical Research
Letters. To appear.

\bibitem{Co1} A.~Connes, \emph{Formule de trace en g\'eom\'etrie
non-commutative et hypoth\`ese de Riemann}, C. R. Acad. Sci. Paris {\bf 323}
(1996), S\'erie I, 1231-1236.

\bibitem{Co2} A.~Connes, \emph{Trace formula in non-commutative Geometry and
the zeros of the Riemann zeta function}, Selecta Math. (N.S.) {\bf 5} (1999) ,
no. 1, 29--106.

\bibitem{Ehm} W.~Ehm, \emph{A family of probability densities related to the
Riemann zeta function}, manuscript (2000), 12 pp.

\bibitem{Gre} U.~Grenander, M.~Rosenblatt, \emph{An extension of a theorem of
    G. Szeg\"o and its application to the study of stochastic processes},
  Trans. Amer. Math. Soc. {\bf 76} (1954), 112-126.

\bibitem{Hof} K.~Hoffman, \emph{Banach spaces of analytic functions\/},
Prentice-Hall, Inc. (1962). (Dover Pub., 1988).


\bibitem{Nev} P.~Nevai, \emph{G\'eza Freud, Orthogonal Polynomials and
Christoffel Functions. A Case Study}, Journ. Approx. Theory {\bf 48} (1986),
3-167.



\bibitem{Nym} B.~Nyman, \emph{On the One-Dimensional Translation Group and
Semi-Group in Certain Function Spaces\/}. Thesis, University of Uppsala, 1950. 55 pp.


\end{thebibliography}
\end{document}